\documentclass[10pt]{amsart}
\usepackage{enumerate,amsmath,amssymb,latexsym,
amsfonts, amsthm, amscd}

\theoremstyle{plain}

\newtheorem{proposition}{Proposition}[section]
\newtheorem{theorem}{Theorem}[section]
\newtheorem{corollary}{Corollary}[section]
\newtheorem{axiom}{AXIOM}

\theoremstyle{definition}

\numberwithin{equation}{section}

\newcommand{\real}{\mathbb R}

\begin{document}
\title[THE COX THEOREM\\UNKNOWNS AND PLAUSIBLE VALUE]
{THE COX THEOREM\\UNKNOWNS AND PLAUSIBLE VALUE}

\author{Maurice J. Dupr\'e}
\address{TULANE UNIVERSTIY}
\email{mdupre@tulane.edu}

\author{Frank J. Tipler}
\address{TULANE UNIVERSITY}
\email{tipler@tulane.edu}


\begin{abstract}
We give a proof of Cox's Theorem on the 
product rule and sum rule
for conditional plausibility without assuming continuity or
differentiablity of plausibility.  Instead, we extend the notion
of plausibility to apply to unknowns giving them plausible values.
\end{abstract}
\maketitle


\section{INTRODUCTION}

Since the work of Laplace \cite{Laplace} in the late 18th century,
there have been many attempts by mathematicians to axiomitize
probability theory.  The most important example in the 20th
century was that of A.N. Kolmogorov \cite{Kolmogorov}, who gave a
very simple measure-theoretic set of axioms that modeled the view
of probability introduced into quantum mechanics by Max Born in
1927.  Remarkably, most physicists, in their non-quantum
applications of probability, have not followed Born or Kolmogorov
but R. T. Cox, who in turn based his approach on Laplace's
original idea that probability theory is a precise mathematical
formulation of plausible reasoning.  These physicists argue that,
while the Kolmogorov axioms are elegant and consistent, they are
much too limited in scope.  In particular, the Kolmogorov axioms
in their original form do not refer to conditional probabilities,
whereas most physics applications of probability theory require
conditional probabilities.  Even though unknown by most
mathematicians who work in probability theory, the Laplace-Cox
approach to probability theory was actually accepted by many
distinguished mathematicians prior to Kolmogorov, for examples,
Augustus de Morgan \cite{Morgan}, Emile Borel \cite{Borel}, Henri
Poincar\'e \cite{Poincare}, and G. P\'olya \cite{Polya}.  For a
discussion of applications of Laplacian probability in the
foundations and interpretation of quantum mechanics see Tipler
\cite{TIPLER}.

Cox's probability theory is not defined by precise axioms, but by
three ``desiderata'': (I) representations of plausibility are to
be given by real numbers; (II) plausibilities are be in
qualitative agreement with common sense; and (III) the
plausibilities are to be ``consistent'', in the sense that anyone
with the same information would assign the same real numbers to
the plausibilities.  Cox (\cite{COX1},\cite{COX}, pg. 16)
purported to show that from these requirements, the plausibilities
satisfied, first the  {\bf PRODUCT RULE}:
$$PL(A \& B|C)=PL(A|B \& C)PL(B|C),$$

\noindent
and the {\bf SUM RULE} $$PL(A|B) + PL(\overline{A}|B) = 1$$

The claim that these two rules follow from the desiderata has come
to be known as {\bf COX'S THEOREM}.  The symbol $PL(A|B)$ means a
conditional plausibility, namely ``the plausibility of $A$ given
that we know $B$.''  The symbol ``$\&$'' represents the logical
``both,'' whereas the bar on top represents logical negation.

We shall give in this paper a rigorous mathematical proof for
Cox's Theorem on the product rule for conditional plausibility of
propositions as used in plausible reasoning, a proof that follows
from precise axioms.  We shall see that our axioms are
mathematically simpler and more intuitive than Cox's desiderata.
In particular, we shall not need to make any continuity or
differentiability assumptions.  It is very important to avoid
assuming continuity if the symbols $A$ and $B$ refer to
propositions --- as they do in Cox's paper and book and as they do
in Jaynes' important book {\it Probability Theory} --- because
propositions are necessarily constructed from a finite number of
symbols, and hence properly belong to the integers and not to the
continuum (as represented, for example, in the G\"odel numbering
scheme in the proof of the G\"odel theorems).  We will not follow Kolmogorov and list a short and ideal set of axioms from which all of probability theory can be derived, but instead give a list of axioms and possible alternatives for several.  All of our alternatives are much less technical and more intuitive than those of such authors as Halpern's for example.  

In order to provide this very simple set of axioms for proving the Cox Theorem and deriving the rules of
probability so as to make them apparent to even a reader without
expert mathematical training, we are led to expand the objective
of plausibility theory to more generally deal with objects we call
unknowns which have plausible values.  The aims of the theory of
plausible reasoning are two-fold. First the aim is to derive the
rules by which logic and common sense constrain our inductive
reasoning in the face of limited information, and second to derive
the rules of probability from simple assumptions, so as to make
them apply to propositions in general. A major motivation in our paper is
to make probability applicable in scientific settings where the
frequency theory of probability is of little or no value, and to
justify a Laplacian or Bayesian approach to
probability(\cite{JAYNES},\cite{SIVIA},\cite{KOOPMAN},\cite{SAVAGE}).
The assumptions need to be well motivated and very simple, and the
proof of the basic rules of probability from these assumptions
should hopefully be trivial. As the counterexample of Halpern
\cite{HALP1} shows, the original assumptions of Cox are
inadequate, and the technical assumptions of Paris \cite{PARIS}
are undesirable and still require an unjustifiable continuity
assumption.  More recently, the work of Hardy (\cite{HARDY},
Theorem 8.1) shows very generally that with sufficient hypotheses
the theorem is true, but the hypotheses on the range of values
could be problematic to verify in practice, even though the
development is important and nontrivial.

In brief, in the standard approach, one assumes a Boolean algebra
$E$ of propositions together with a real valued function $PL(A|B)$
defined for $(A,B) \in E \times E_0,$ where $E_0=E \setminus
\{0\}$ and which we think  of $PL(A|B)$ as assigning a numerical
level of {\bf PLAUSIBLILITY} to proposition $A$ given that we
accept proposition $B$ as true. We further assume that $PL(A|B)$
is a monotonic function of the plausibility of $A$ when $B$ is
assumed true. Consequently, we are allowed to modify $PL$ by
composing with a monotonic function if necessary to produce a
useful rule. The first result of the standard approach is that
allowing such modifications we can produce the product rule and
the sum rule. In the original proof of the product rule, which is
essentially the rule for conditional probability, R.T. Cox (see
\cite{COX}, page 12) assumed merely that the plausibility $PL(A \&
B|C)$ was a numerical function of the plausibilities $PL(A|B \&
C)$ and $PL(B|C)$ through some real valued function $F$ of two
variables. Motivating this assumption requires examination of a
host of special cases \cite{TRIBUS1} for the different
possibilities of what $PL(A \& B|C)$ could depend on among the
four numbers
$$PL(A|B \& C), ~PL(A|C),~PL(B|A \& C), ~PL(B|C),$$ an examination
rendered unnecessary in the approach we will introduce here. Then
evaluating $PL(A \& B \& C|D)$ in the two possible ways available
and applying associativity of the conjunction of propositions
almost leads to the conclusion that the function $F$ is an
associative multiplication on the set of real numbers forming the
range of $PL.$ This last step taken by Cox was a logical mistake
as the counterexample of Halpern \cite{HALP1} shows, this
conclusion is not justified as $E$ may be finite, and even if it
were true, there would in general be no useful information coming
from this fact. But Cox assumed that the function $F$ should be of
a universal character and therefore must be defined on the whole
plane. Cox thus assumed $F$ is an associative multiplication on an
interval of real numbers. Assuming the function to be
differentiable leads to the assumed multiplication being in fact
ordinary multiplication. However, the assumption that the function
is differentiable was never justified by either Cox or Jaynes,
except by hand waving. Moreover, the domain of the function may in
reality only be a finite set of real numbers, and so extreme
effort has gone into trying to add on very technical assumptions
\cite{PARIS}, \cite{HALP1}, which in effect produce sufficient
density of the domain to claim that continuity gives associativity
which together with strict monotonicity (a requirement from
"agreement with common sense") suffices to show that the
associative multiplication is just ordinary multiplication.

However, it has been well known for many many years by experts in
the theory of topological semigroups, what possible continuous
multiplications are available on an interval of real numbers.
Numerous textbooks in topological semigroup theory address this
very issue.  As shown for example in the seminal work by K.H.
Hofmann and P.S. Mostert \cite{KHH}, the possibilities are
infinite. However, if we assume the strict monotinicity which
seems consistent with common sense and which rules out idempotents
other than a zero and a unit to form the boundary of the interval,
then the only continuous associative multiplication is isomorphic
to the unit interval under ordinary multiplication. Thus a
suitable function of $PL$ would then satisfy the multiplication
rule, a result which properly belongs to the theory of topological
semigroups. Several authors have dealt with counterexamples
\cite{DUB}, \cite{HALP1}, \cite{KRAFT} and proofs \cite{PARIS}, in
effect reproving results of topological semigroup theory, and the
complete proof of Cox's Theorem, even with the assumption of
continuity, is not simple. In the case of Hardy \cite{HARDY}, we
have a fairly complete theory of scales which in effect provide
alternate density type assumptions on the set of values of the
plausibility function (\cite{HARDY}, Theorem 8.1). The scales are
themselves lattices of special type which under the proper
technical assumptions are shown isomorphic to the unit interval.
This approach is very general and in spirit similar to the
(noncommutative operator algebra) case treated by Loomis
\cite{LOOMIS}. Moreover, these technical assumptions such as
continuity or divisibility are just as problematic as the
assumption of differentiability. That is, they are certainly
reasonable, and not as strong as assuming differentiability, but
in the end, they are still strong and highly technical,
non-intuitive assumptions. Once the Cox Theorem is proved, the
modified function $PL,$ under another common sense assumption,
namely that $PL(notA|B)$ depends only on $PL(A|B)$, can be shown
to have (at least a power, depending on which axioms are used)
which obeys the laws of probability. Several authors have dealt
with the problem of associativity of the universal function,
required for the Cox Theorem, as it seems to be essential to the
argument given by Cox, and it is essentially a result of
topological semigroup theory which is being applied by all these
authors. However, we will see that in our approach, questions of
continuity or associativity become completely irrelevant to the
argument.

Our simpler approach takes a closer look at what scientists are
really trying to do.  The main aim of scientists, engineers and
technical workers is arriving at values for numerical quantities
on the basis of limited information.  Thus, instead of restricting
attention to a set of propositions, we are led instead to consider
a set of more general objects we shall call {\bf UNKNOWNS}.  We
purposefully do not use the term ``random variable'' here, as it
is a much too restrictive a notion, and carries with it all the
baggage of the Kolmogorov approach to probability theory, but a
random variable is an example of an unknown. In case of
propositions, since all members of a Boolean algebra are
idempotent, and as the only idempotent numbers are 0 and 1, we are
naturally lead to create or define an unknown number, $I_A,$ for
each proposition $A$ called its {\bf INDICATOR}. Our object now is
to assign a {\bf PLAUSIBLE VALUE} denoted $PV(X|A)$ to the unknown
$X$ given the information in proposition $A.$  As for plausibility
of statements, we then simply define $PL(A|B)=PV(I_A|B).$ The
result is we find a very simple and natural theory of plausible
value for unknowns which contains the theory of plausibility of
propositions and which requires no assumptions at all in the form
of differentiablility or even of continuity for its rules.  The
rules are simply dictated by simple common sense consistency with
logic. The main idea turns out to be exceedingly simple and really
only depends on some simple properties of retraction mappings on
sets. What comes out of these considerations is that the rules are
really uniquely determined, in a very strong sense, merely by the
assumption that some form of rule exists. In short, existence
implies a strong form of uniqueness. We begin with simple
considerations of retraction mappings on sets, and then when we
get to the setting of unknowns, we see right away that the $PV$
must be a retraction of the unknowns onto the knowns.  Thus, the
assumption of the existence of rules of dependency of certain
general forms can be completely determined by what happens to the
known quantities under the general forms of the rules.  In
particular, if we examine what this approach does for the
plausibility theory, we note that a natural logical axiom of
rescaling of plausible value under changes of units causes the
universal function of the Cox theorem proof to be homogeneous in
its first variable. This axiom for plausibility means that
plausibility should really be a geometric quantity which is
independent of the choice of maximum and minimum.  That is, we
should think of the plausibility of a statement as being specified
by a point on a line segment where one endpoint is the
plausibility of a known true statement and the other endpoint is
the plausibility of a known false statement. That geometric
picture is independent of the  numerical scale chosen for the
segment, and a realistic plausibility theory should contain that
property.  That is, if someone asks you what is the plausibility
of statement $A$ given statement $B$ is true on a scale of $a$ to
$b$, you should be able to express the plausibility on that scale
demanded no matter what scale you had originally chosen to express
plausibility.  What this means is that if we define
$O(A|B)=PL(A|B)/PL(notA|B),$ usually called the odds of $A$ given
$B,$ then $O(A|B)$ is completely scale invariant.  Homogeneity of
the universal function of the Cox theorem gets around the
counterexample of Halpern \cite{HALP1}. In fact by Halpern's
theorem 3.1 and lemma following, if $F$ satisfying the conclusion
of his theorem is homogeneous in the first variable, then we find
immediately that $F(x,y)=xy$ as an immediate consequence of his
theorem 3.1, so $F$ is associative, a contradiction of his
following lemma.  His construction technique is to take a finite
set of 12 members and by using two slightly different probability
distributions, join them in an unnatural way to produce a
plausibility theory which satisfies the assumptions of Cox but for
which the universal function $F$ cannot possibly be associative
because of the way the two probability distributions are joined to
produce the plausibilities.  Of course, we see immediately now,
that Halpern's counterexample violates the natural rescalability
that plausibility should have, that is, his function $F$ cannot be
homogeneous in its first variable, so his counterexample fails to
be a counterexample in any system of plausibility theory in which
plausibilities have a natural scale invariant meaning.

\section{SIMPLE RETRACTION PRINCIPLES}

One of the first things a mathematics student learns is that if
$f$ and $g$ are functions on the set $T,$ if $g$ has range $S$ so
that $g(T)=S,$ then there is at most one function $h$ with domain
$S$ satisfying $f=hg.$  In short, for such $h$ to exist, clearly
$f(t)$ as a function of $t \in T$ must only depend on the value
$g(t),$ or in other words, if $g(t_1)=g(t_2),$ then
$f(t_1)=f(t_2).$  If we assume this condition is satisfied, then
using the axiom of choice if necessary, we can form a {\bf
SECTION} of $g,$ namely a function $s$ from $S$ to $T$ with the
property that $gs=id_S,$ the identity function on $S.$  We get $h$
on setting $h=fs.$  For then, $hg=fsg,$ but $gsg=g$ implies
$fsg=f,$ by the assumed condition.  In a sense here, we can say
existence implies uniqueness, but the function $h$ we find does
not have a simple dependence on $f$ for its construction. We may
have to use the axiom of choice. We will see that the dependence
of $h$ on $g$ is quite explicit if $g$ is a retraction onto a
subset of $T.$

To begin, recall that if $T$ is any set,  $R \subset T$ is any
subset of $T,$ then, a {\bf RETRACTION} $P$ of $T$ onto $R$ is a self
mapping of $T$ such that its image is $R$ and $P(x)=x$ for each $x
\in R.$ We shall also find it useful to recall the idea of a {\bf RESTRICTION}
of a function: if $f$ is a function defined on $T,$ then we denote
by $f|R$ its restriction to the subset $R,$ that is the same rule,
but with domain restricted to be $R.$

\begin{proposition}\label{retractfactorize}
Suppose that $P$ is a retraction of the set $T$ onto the subset
$R$ and that $f$ is a function from $T$ to set $S.$  If $f(t)$ for
$t \in T$ only depends on the value $P(t),$ then there is a unique
function $h$ defined on $R$ with $f=hP,$ and in fact $h=f|R,$ the
restriction of $f$ to $R.$
\end{proposition}

\begin{proof}
The hypothesis that $f(t)$ only depends on $P(t)$ guarantees the
existence of $h.$ But now, for $r \in R,$ we have $P(r)=r$ as $P$
is a retraction onto $R,$ and hence $f(r)=h(P(r))=h(r),$ so
$h=f|R.$
\end{proof}

\begin{corollary}\label{retractcor}
Suppose that $P_k$ is a retraction of the set $T_k$ onto the
subset $R_k,$ for $k=1,2,3.$ Suppose $m$ is a mapping from $T_1
\times T_2$ into $T_3$ with $m(R_1 \times R_2) \subset R_3,$ and
denote this mapping by juxtaposition, $m(x,y)=xy.$  Then:

(1) if $f$ is a function from $T_1$ to $T_2$ with $f(R_1) \subset
R_2,$ and if $P_2(f(t))$ depends only on $P_1(t),$ then

\begin{equation}
P_2(f(t))=f(P_1(t)), ~~~t \in T.
\end{equation}

(2) if $P_3(t_1t_2)$ depends only on $(P_1(t_1),P_2(t_2)),$ then

\begin{equation}
P_3(t_1t_2)=P_1(t_1)P_2(t_2)),~~~(t_1,t_2) \in T_1 \times T_2;
\end{equation}

(3) if in (2) we have a fixed $e \in T_2$ and if we instead assume
that $P_3(t_1e)$ depends only on $P_1(t_1),$  then

\begin{equation}\label{coxgen1}
P_3(t_1e)=P_3([P_1(t_1)]e),~~~t_1 \in T_1;
\end{equation}

(4) if for (3)  in  addition we assume $e$ has the property that
$P_3(re)=rP_2(e),$ for all $r \in R_1,$ then

\begin{equation}\label{coxgen2}
P_3(t_1e))=P_1(t_1)P_2(e),~~~t_1 \in T_1.
\end{equation}

\end{corollary}

\medskip

\begin{proof}
The hypothesis in (1) guarantees a function $h$ defined on $R_2$
with the property that $P_2f=hP_1.$  But now the proposition tells
us that $h=P_2f|R_1,$ but $P_2f|R_1=f|R_1,$ because $f(R_1)
\subset R_2$ and $P_2$ is a retraction onto $R_2.$  The hypothesis
in (2) guarantees that $P_1 \times P_2$ is a retraction of $T_1
\times T_2$ onto $R_1 \times R_2,$ and hence using (1) with $f=m$
completes the proof for (2).  In case of (3), with $e \in T_2$
fixed, we have a unique function $h_e$ from $R_1$ to $R_3$ such
that $P_3(t_1e)=h_e(P_1(t_1)),$ for all $t_1 \in T_1.$  But then,
taking $r \in R_1,$ we have $P_1(r)=r,$ so

$$h_e(r)=h_e(P_1(r))=P_3(re),$$ and [\ref{coxgen1}] follows
immediately. Now, (4) is clear from (3).

\end{proof}

\medskip

In particular, if we take $P_1=P_2=P$ in (1) of the corollary,
then we see that $P(f(t))$ depends only on $P(t)$ exactly when
$Pf=fP,$ a {\bf GENERAL COMMUTATION RULE}.  In case of (2), we
have a {\bf GENERAL COMBINATION RULE:}  if $P_3(xy)$ depends only
on $(P_1(x),P_2(y),$ then $P_3(xy)=P_1(x)P_2(y).$ On the other
hand, if we take the case where $T_1=R_1,$ so $P_1$ is simply the
identity on $R_1,$ then when $P_3(ry)$ depends only on
$(r,P_2(y))$ for $r \in R_1$ we conclude from (2) that
$P_3(ry)=rP_2(y),$ a form of {\bf GENERAL HOMOGENEITY}.  We can
also conclude this for fixed $r$ in $R_1$ using (1). That is, we
take $f$ above to be left multiplication by $r \in R_1.$  We can
note that (4) above is a very general form of the product rule part of Cox's Theorem.  In
particular, we note that the question of any form of associativity
never enters the proof of (4).

\section{UNKNOWNS AND PROPOSITIONS}

Scientists, engineers, and technical workers deal with a world of
numbers, and other mathematical entities many of which are not
completely known. In many situations, when the description of a
particular quantity's numerical value tells us only that a well
defined value exists without telling us what it actually is, we
must proceed with a most plausible value based on the information
at hand which may be incomplete, and which may not be certain. The
information generally appears as a proposition which in fact is
either true or false, and once accepted is assumed true for
purpose of evaluating the unknown quantity as well as we can. Such
quantities are actually more than simple real numbers, as their
descriptive information is part of their structure and does not
generally give us enough information to determine a certain value.
Thus, we can consider them to be objects in some set containing the
set of all real numbers and that there is some real valued
function on that set which gives each object a value and that this
function is unknown to us. We wish to analyze how the requirement
of logical consistency constrains the procedure for arriving at
plausible values for these objects or unknown quantities when
limited information is available.  Even if we are just guessing,
their should be certain simple logical constraints. As Cox
\cite{COX} has shown, if we try to apply plausibility with no
information, we arrive at absurd results, so our prior information
must give us some information about an unknown of interest. More
generally, scientists and engineers often deal with mathematical
structures beyond the real number system and the same
considerations apply. When a physicist speaks of the state
$\varPsi$ of a classical bounded quantum mechanical system, he
generally means that $\varPsi \in H,$ where $H$ is some Hilbert
space, but before he applies the rules of quantum mechanics, he
really does not know what $\varPsi$ is. In fact, he may not even
know what $H$ is. In fact, he may not know enough about the actual
physical system for the rules of quantum mechanics to determine
what $\varPsi$ is. He assumes by the axioms of quantum mechanics
that the physical system under consideration determines a unique
state, but the information and measurements he actually has for
the system may not be enough to actually determine $\varPsi.$ For
instance, $\varPsi$ could be the state of a black cat in a closed
box which we cannot see inside, but which we can hear meowing. We
could therefore properly think of $\varPsi$ as a symbol for an
unknown unit vector in $H,$ and we could try, based on $C,$ the
proposition stating the measurements we have made and our
knowledge of quantum mechanics, to arrive at a plausible value
$PV(\varPsi|C) \in H.$  The same type of consideration applies to
any unknown member of any set based mathematical structure.
Information can appear in the form of differential equations which
must be satisfied as well as experience we have in dealing with
similar problems in our past-everything we know can be brought to
bear on the choice of a plausible value.  When the mathematical
structure has rules of combination such as vector addition,
semigroup multiplication, actions of one system on another, and so
forth, clearly these same operations should apply to the unknowns.
Thus, if $X$ is an unknown number and $\varPsi$ is an unknown
vector in $H,$ then $X \varPsi$ is another unknown vector in $H.$
If we are interested in the unknown $\varPsi$ in $H$ and the
unknown $\varPhi$ in $H,$ then we possibly we could end up needing
to consider $\varPsi + \varPhi.$  Certainly if we have information
about each of the summands, then we know something about the sum.
Thus it is reasonable to assume that whatever unknowns we are
interested in dealing with algebraically form the same kind of
system as the system they "live in". For instance, we could think
of the Hilbert space $H$ as being an unknown member of a small
category of Hilbert spaces if it is also unknown.

To begin, let us be precise about our set up and then consider examples of what we mean by an
{\bf UNKNOWN}.  Suppose that $S$ is any set.
Suppose that $B$ is a proposition which describes a member $X$ of
$S$ sufficiently well so that $B$ implies such a member exists even
though $B$ might not state which member of $S$ it is, then $X$ is
an unknown member of $S.$  In particular, if $s \in S,$ then we
regard $s$ as known, that is, a known unknown.  Thus, if we are
interested in a set $T$ of unknown members of $S,$ then we usually
assume that $S \subset T.$ That is to say, we should think of the
unknowns in $S$ as having additional structure by virtue of their
descriptions, and we regard the known members of $S$ as contained
in the unknowns.  To proceed formally, then we will simply assume
that $S \subset T$ are sets and we are regarding $T$ as the set of
unknowns of $S$ in which we are interested.  Of course, as each $X
\in T$ is an unknown member of $S,$ it must have a value $AV(X),$
called the {\bf ACTUAL VALUE} of $X,$ but we are in general not
aware of what this is. That is, we have limited information about
it. Of course, $AV(s)=s$ for each $s \in S,$ that is we assume the
members of $S$ are trivially known. The plausible value function
$PV$ is mathematically an $S$-valued function defined on $T \times
E_0,$ where $E$ is a Boolean algebra of propositions and $E_0$
denotes the non-zero members of $E.$  We denote by $PV(X|A)$ the
value of this function on the pair $(X,A) \in T \times E_0.$

We must make some basic assumptions on how unknown quantities get
plausible values. Now, the most basic assumption that can be made
which is absolutely obvious from the standpoint of logical
consistency is that if our information tells us exactly what value
an unknown has, then the plausible value of that unknown given
that information must be that value the information is telling us.
So we formulate this as our {\bf FIRST AXIOM OF PLAUSIBLE VALUE}.

\medskip

\begin{axiom}\label{axpv}
If $T$ is a set of unknown members of the set $S,$ where $S$ is
any set, we assume that $S \subset T$ and $AV$ is a retraction of
$T$ onto $S$. If $X \in T,$ if $s \in S$ and the proposition $A
\in E_0$ implies that $AV(X)=s,$ then $PV(X|A)=s.$
\end{axiom}

\medskip

Notice by Axiom \ref{axpv} of plausible value, that $PV(\_|A)$ for
fixed proposition $A \in E_0$ defines a retraction of $T$ onto
$S,$ if $T$ is a set of unknown members of $S.$ This is because if
$s \in S,$ then $A$ trivially implies $AV(s)=s$ so by axiom 1 we
have $PV(s|A)=s.$

Our next axiom also makes good common sense from the standpoint of
logic. If our information is telling us that two unknowns have the
same value, even if we do not know that value, we must choose the
same plausible value for both in order to maintain logical
consistency.

\medskip

\begin{axiom}\label{axeq}
If $X,Y \in T$ are unknown members of the set $S$ and if the
proposition $A \in E_0$ implies that $AV(X)=AV(Y),$ then
$PV(X|A)=PV(Y|A).$
\end{axiom}

\medskip

And now for the examples.  Consider a set $S$ and any set
$D$ and form the set $T$ of $S-$valued functions on $D,$ so
$T=S^D.$ We regard $S \subset S^D=T$ by identifying each member of
$S$ with a constant function on $D.$ Let $d \in D$ and define
$AV(X)=X(d)$ for each $X \in T.$  Of course, taking $PV=AV$
independent of the $E_0$ variable satisfies the axioms showing
consistency.

In particular, consider unknown real numbers. We regard
an {\bf UNKNOWN (NUMBER)} as any defined numerical quantity $X$
whose definition tells us it has an exact value but whose
definition does not necessarily tell us what that value is.
Suppose we have some assumed information in the form of a
proposition $C$ which influences our idea of what its value might
be.  For example, $X$ could be Beethoven's weight in pounds at
noon on his fifth birthday. We can take $C$ to be a proposition
which states our knowledge of typical weights of five year old
children. Clearly 1000 is not a reasonable guess as to what $X$
is, but 45 might not be to far off. As another example, we can
take $Y$ to be the current outside temperature in degrees Celsius.
If $C$ is the statement of all of our previous knowledge of
weather, our experience of the outside air temperature the last
time we were outside, as well as what we see by looking out our
office window, then we may be able to get a pretty good plausible
value of the outside temperature.  If we are outside we can
probably do even better.  Now, our plausible value may be only a
guess, and there may be many choices, but we want to imagine that
there is some set $E$ of propositions that we will consider and
some set $T$ of unknowns that we are interested in, and that for
these we choose $PV(X|C)$ for each $C$ in $E_0$ and each $X$ in
$T.$ Now, again, we want to develop the properties of $PV$ based
on the idea that as a function on $T \times E_0$ to ${\real},$ it
must have certain properties to conform to common sense logical
consistency.

We can notice that if $X$ and $Y$ are unknown numbers, then we can
clearly form $X+Y$ and $XY.$ For instance, $X$ and $Y$ could be
the unknowns in the two preceding examples involving weight and
temperature. If we have some information about $X$ and $Y,$ then
we have information about their sum and product as well. The
unknowns have no units in and of themselves, the units are
contained in their descriptive information which gives them a
numerical value, so any unknown numbers can always be added and
multiplied.  Since it is reasonable to assume that if we are
interested in a pair of unknown numbers we might also need to deal
with their sum and product, we assume then that $T$ is closed
under the operations of addition and multiplication, making it a
{\bf RING}.  This is mainly a convenience, and we should point out
that for our proof of the Cox Theorem, we only need to assume
closure under multiplication of unknowns by indicators, which we
proceed to define next. We assume that if $A$ and $B$ belong to
$E,$ then so do $A \& B,$ the negation of $A,$ denoted $notA,$ and
$A ~or~ B$ and that $E$ is nonempty, so it is a {\bf BOOLEAN
ALGEBRA} of propositions. If $C$ is a proposition, then we can use
it to define an unknown $I_C$ which has the value 1 if $C$ is true
and the value 0 if $C$ is false, and which we call the {\bf
INDICATOR UNKNOWN} of $C.$ Notice the truth value of a proposition
is entirely contained in its indicator unknown, so interest in
whether or not a particular proposition is true is equivalent to
interest in the value of its indicator unknown. Consequently, we
assume that $T$ contains all indicators of propositions in $E.$ As
with general sets, we will regard the real scalar field, ${\real}$
as special unknowns which are known values under any information,
($C$ implies $AV(r)=r$ for every number $r$), so we assume that
$T$ contains ${\real},$ the field of real numbers and therefore in
particular, $T$ is an {\bf ALGEBRA} over ${\real}.$   As far as
the Boolean algebra $E$ is concerned, we can note that in general,
by Stone's Theorem \cite{HALMOS}, we can embed $E$ as a Boolean
algebra of idempotents in the algebra $C_E$ of continuous real
valued functions on the Stone space of $E.$  We can therefore
regard the algebra $T$ as an algebra over $C_E$ as a way of more
concretely thinking of the way indicators act on unknowns.
Similarly, if $W$ is a vector space and $T$ is a vector space of
unknown members of $W,$ then we can regard the action of
indicators on $T$ as coming from a $C_E-$module structure on $T.$
Thus, if $K$ is any commutative algebra over ${\real},$ then we
can take any $K-$module $T$ with ${\real}-$submodule $W,$ a
retraction $AV$ of $T$ onto $W,$ and for each idempotent $A$ in $K
\setminus 0$ choose a retraction $PV(\_|A)$ of $T$ onto $W,$ to
produce a mathematical model of the setup for unknown vectors in
$W.$ Since these retractions can be chosen to be linear, we see
that there exist many such setups.

We summarize these comments as our next axiom.

\medskip

\begin{axiom}\label{axrealunknown}
We assume a set $T$ of real unknowns  is a commutative algebra
with identity over the field of real numbers, ${\real},$ and that
it contains the indicator unknowns of all propositions in the
Boolean algebra of propositions $E,$ that is we assume that the
set of indicators of members of $E$ is a Boolean algebra of
idempotents in $T.$
\end{axiom}

\medskip

We want to put order axioms on our plausible numerical values so
that plausible numerical values are logically consistent with
common sense. In particular, we will take as our next axiom:
\medskip

\begin{axiom}\label{axorder}
If $X$ and $Y$ are in $T$ and if $C$ is in $E_0,$ and if $C$
implies that $AV(X) \leq AV(Y),$ then $PV(X|C) \leq PV(Y|C).$
\end{axiom}

\medskip

This axiom merely says that we must choose the ordering of
plausible values so as not to contradict the order information we
have about the underlying numerical unknowns. As an immediate
consequence of this axiom, we have that if $C$ implies that
$AV(X)=AV(Y),$ then $PV(X|C)=PV(Y|C).$  This is simply because for
real numbers, $=$ is the same as $\leq \& \geq.$ Thus, we see that
Axiom \ref{axeq} in the case where $S={\real},$  is a consequence
of Axiom \ref{axorder} for the case where $S={\real}.$ In
particular, as a consequence of Axiom \ref{axpv} , if $r$ is any
real number, then since $C$ trivially implies $AV(r)=r,$ it
follows that $PV(r|C)=r.$   Thus for fixed $C,$ the plausible
value $PV(X|C)$ viewed as a function of $X$ in $T$ is in fact a
retraction of $T$ onto ${\real} \subset T.$ Now an immediate
consequence of Axioms \ref{axpv} and \ref{axorder} is that if $a$
and $b$ are real numbers and $C$ implies that $a \leq X \leq b,$
then $$a \leq PV(X|C) \leq b.$$ If $A,C$ are in $E,$ then $0 \leq
I_A \leq 1,$ so by Axioms \ref{axpv} and \ref{axorder} we can
immediately conclude that
$$0 \leq PV(I_A|C) \leq 1.$$
In view of the preceding inequality, we define the {\bf
PLAUSIBILITY} of $A$ given $C,$ denoted $PL(A|C),$ by
$$PL(A|C)=PV(I_A|C).$$

Now, it is certainly reasonable that if $X$ is in $T$ and we have
determined $PV(X|C)$ and if $r$ is any real number then we should
be able to determine $PV(rX|C)$ from $r$ and the purely numerical
value $PV(X|C).$ For instance, we should be able to change units
and do unit conversions directly on the plausible values (if you
think the plausible value for the outside temperature is 20
degrees Celsius, then you should think it is 68 degrees
Fahrenheit).  At least we should be able to rescale plausible
values under unit changes, even if we do not accept changes of
zero point as in temperature conversion. This leads to our next
axiom:
\medskip

\begin{axiom}\label{axhomog}
If $r$ is any real number, if $C$ is any proposition in $E_0,$ and
if $X$ and $Y$ are unknowns in $T,$ and if $PV(X|C)=PV(Y|C),$ then
$PV(rX|C)=PV(rY|C).$  In other words, we assume that $PV(rX|C)$
for fixed $r \in {\real}$ depends only on $PV(X|C).$
\end{axiom}
\medskip
Thus, by (1) of corollary [\ref{retractcor}], we have homogeneity
of plausible value:
\begin{equation}\label{realhomog}
PV(rX|C)=rPV(X|C)
\end{equation}

\medskip

Finally, we consider the axiom that leads to our form of Cox's
Theorem which we shall call the {\bf COX AXIOM:}
\medskip

\begin{axiom}\label{axcox}
If $A,C$ are fixed in $E,$ if $X_1,X_2$ are in $T,$ if $PV(X_1|A
\& C)=PV(X_2|A \& C),$ then $PV(X_1 I_{A}|C)=PV(X_2 I_{A}|C).$
That is, we assume that as a function of $X,$ the plausible value
$PV(XI_A|C)$ depends only on $PV(X|A \& C).$
\end{axiom}
\medskip

To motivate this axiom, notice that if $A$ is false, then
$XI_A=0,$ whereas if $A$ is true, then we are evaluating the
plausible value of $X$ with both $A$ and $B$ being true, which
should somehow depend only on $PV(X|A \& C).$  Notice the
asymmetry here, which prevents any consideration of the multitude
of possibilities in plausibility theory \cite{TRIBUS1}. We cannot
put $X$ in the position of the given information, the first
variable of $PV$ can only be an unknown and the second variable
can only be a statement. Moreover, $PV(XI_A|C)$ cannot depend on
the numerical value of $PV(X|C)$ because we could generally have
unknowns $X$ and $Y$ with $PV(X|C) \neq PV(Y|C)$ but with $A$
implying that $X$ and $Y$ are equal, in which case we clearly must
have that $PV(XI_A|C)=PV(YI_A|C).$ This leads directly to our form
of the product rule of Cox's Theorem.

\medskip

\begin{theorem}\label{coxthmpv}
If $X$ is any unknown number in $T$ and if $A,C$ are any
propositions in $E,$ with $A \& C \in E_0,$ then

\begin{equation}\label{coxthmpvrule}
PV(XI_A|C)=PV(X|A \& C)PV(I_A|C).
\end{equation}
\end{theorem}

\medskip

\begin{proof}
This is an immediate consequence of (4) in corollary
(\ref{retractcor}) and the previous axioms, where we take
$P_3=PV(\_|C),~~P_1=P(\_| A \& C),$ and $P_2=PV(\_|C).$
\end{proof}

\begin{corollary}\label{coxthmpl}
If $A,B,C$ belong to $E,$ with $B \& C \in E_0,$ then

\begin{equation}\label{coxthmplrule}
PL(A \& B|C)=PL(A|B \& C)PL(B|C).
\end{equation}

\end{corollary}

\noindent
which is the standard product rule of Cox's Theorem.

We need to point out here, that our approach to the Cox theorem
(\ref{coxthmpl}) has eliminated the problems which allow the
counterexample of Halpern \cite{HALP1}.   We do not need to have
an associative multiplication on the real line or an interval, we
do not need to assume any continuity or differentiability or
divisibility, we do not need to assume that our Boolean algebra of
propositions has sufficiently many plausible values to have dense
range in an interval of numbers. We do not even need to assume a
function of two real variables as Cox does, we merely assume that
for fixed $A \in E$ that the plausible value of $XI_A$ as it
depends on $X$ is somehow only depending on the plausible value
assigned to $X,$ a considerably weakened assumption. In fact, we
could have the hypothesis only for a particular $A$ and the result
then applies to that particular $A.$ That is, by
(\ref{retractcor}), we see that we do not even need to assume this
for all $A \in E$ at once, it is enough to assume it for a single
$A \in E$ and to assume the homogeneity of that single indicator.
In effect, by passing to unknowns and using indicator functions,
the rescalability of plausiblilities (encoded in the homogeneity
of plausible value) causes the range to be the whole real line and
consequently, the universal function assumed by Cox will here have
to have domain ${\real} \times Im(PL).$  It is thus the
homogeneity, not the additivity assumed by some authors
\cite{KOOPMAN}, which is the crucial ingredient which gives the
result. Also, we have assumed our set of unknowns forms an algebra
with identity over the reals as it seems most natural, but we
really only used the fact that we have a set of unknowns that is
closed under scalar multiplication and contains all indicators
from $E$ and the real numbers themselves as a subset. We can also
point out that if we drop the axiom of homogeneity
(\ref{axhomog}), then by (3) of the corollary (\ref{retractcor}),
we would still obtain a weakened form of the Cox Theorem as a
consequence of the other axioms.

An alternative to the Cox Axiom, due to Savage \cite{SAVAGE}, in
case of plausibility, is the {\bf SURE THING AXIOM:}

\medskip

\begin{axiom}\label{axsurething}
If $X,Y \in T$ and $A ,B,~A \& B, B \setminus A \in E_0,$ and if
both $PV(X|A \& B)=PV(Y|A \& B)$ and $PV(X|B \setminus A)=PV(Y|B
\setminus A),$ then $PV(X|B)=PV(Y|B).$
\end{axiom}

\medskip
If we form $Y=PV(X|A \& B)I_A \in T,$ then the sure thing axiom
implies that $PV(XI_A|B)=PV(Y|B),$ since $PV(XI_A|A \& B)=PV(X|A
\&B)$ by axiom \ref{axeq}, and the product rule of Cox's Theorem is then an
immediate consequence of this equality and homogeneity from axiom
(\ref{axhomog}).

So far, nothing has been said about additivity of $PV.$  Of
course, (2) of corollary (\ref{retractcor}) gives additivity if we
assume there is an appropriate general dependence.

\begin{proposition}
Suppose that $S$ is a set with binary operation, $+,$ and $T$ is a
set of unknowns of $S,$ which is closed under $+.,$ and with $S
\subset T.$ If we assume that $PV(X+Y|A)$ for all unknowns $X$ and
$Y$ in $T$ depends only on the values $PV(X|A)$ and $PV(Y|A),$
then

\begin{equation}\label{gensum}
PV(X+Y|A)=PV(X|A)+PV(Y|A).
\end{equation}
\end{proposition}

\begin{proof}  This is an immediate consequence of (2) in
corollary [\ref{retractcor}] on taking $P_1=P_2=P_3=PV(\_|A).$
\end{proof}

What we see here is that the additivity of plausible value in the
most general sense possible would be a consequence of the basic
logical consistency of meaning together with the mere assumption
that some form of law of combination exists.  Thus, the same would
apply if we were considering plausible values of unknown vectors
in vector spaces-if we assume the plausible value of the sum
somehow depends on the plausible value of the summands, then the
only possible rule is the standard sum rule. Such general
additivity laws are usually easy to motivate with examples, or in
the case of $S={\real}$ by thinking in terms of money, but in the
end, whatever the motivation, it includes the motivation that an
actual rule exists, and that is already enough. Moreover, the
final arbiter on such an assumption has to be whether experience
with its use leads to reasonable results.  For, notice that if the
operation is taken to be ordinary multiplication with $S={\real},$
the same argument applies but then the rule is not generally true
even for ordinary expectations in ordinary probability theory,
which means that for general expectations in probability theory
there can be no general rule for getting the expected value of a
product from the expected values of the factors.  We see from
proposition (\ref{gensum}) that if we axiomatically assume there
is some form of rule giving the plausible value of a sum in terms
of the plausible values of the individual summands, then the only
possible rule is the ordinary sum rule.  But, before going that
far, let us reconsider the temperature example.

Suppose that $X$ is the outside temperature in degrees Celsius. If
our information leads to a best guess of $c$ as the most plausible
value, then consistency requires that in degrees Fahrenheit the
plausible value is $32+(9/5)c.$  This includes a change in zero
point. Thus, consistency with the most general changes of units
for any  unknown leads to the next axiom:
\medskip

\begin{axiom}\label{axrescale}
If $a,b$ belong to ${\real},$ if $X$ belongs to
$T,$ and if $C$ belongs to $E,$ then
\begin{equation}\label{unitchange}
PV(aX+b|C)=aPV(X|C)+b.
\end{equation}
\end{axiom}
\medskip

Notice that this axiom implies axioms \ref{axpv} and \ref{axhomog}
and includes a limited form of additivity.  Thus, in particular,
axiom \ref{axcox} and this axiom imply the sum rule of Cox's theorem.  However,
this last axiom allows us to immediately arrive at the properties
of plausibility for statements.  Because we have

$$I_{notA}=1-I_A$$
$$I_{A \& B}=I_A I_B$$
and therefore by de Morgan's Law

$$I_{A or B}=I_A + I_B -I_A I_B.$$
So,
$$PV(I_{notA})=1-PV(I_A),$$
and it is well known \cite{JAYNES} that the sum rule of Cox's theorem and the
preceding complementation property imply by deMorgan's Law that
$PV(\_|C)$ is additive on indicators of exclusive propositions. We
thus arrive at the usual rules of probability on defining the
probability, $P(A|B),$ of $A$ given $B$ by $P(A|B)=PV(I_A|B).$ To
obtain the general additivity of plausible value, we now only need
to assume the following simpler axiom.

\medskip

\begin{axiom}\label{axadd}
For $T$ an algebra of real unknown numbers, for each fixed $Y \in
T$ and $A \in E_0,$ the plausible value $PV(X+Y|A)$ depends only
on $PV(X|A).$
\end{axiom}

\medskip

\begin{proposition}\label{pvaddreal}
If $T$ is an algebra of unknown numbers and $X,Y \in T$ with $A
\in E_0,$ then assuming axioms \ref{axrescale} and \ref{axadd},

\begin{equation}\label{realadd}
PV(X+Y|A)=PV(X|A)+PV(Y|A).
\end{equation}

\end{proposition}

\medskip

\begin{proof}
Fix $Y \in T$ and $A \in E_0.$  Now, by assumption, on considering
$PV(X+Y|A)$ as a function of $X$ alone, the Axiom \ref{axadd}
guarantees a function $f_{(A,Y)}$ satisfying
$f_{(A,Y)}(PV(X|A))=PV(X+Y|A),$ for every $X \in T.$  If we take
the special case of $X=r \in {\real},$ then, as $PV(r|A)=r,$ and
as by Axiom \ref{axrescale} we have $PV(r+Y|A)=r+PV(Y|A),$ it
follows that
$$f_{(A,Y)}(r)=PV(r+Y|A)=r+PV(Y|A),$$
for every real number $r,$ and this gives the result.
\end{proof}

\medskip

We are of the opinion that the most economical approach to
probability theory is to take as axioms, \ref{axorder},
\ref{axcox}, and \ref{axrescale}, as these three axioms easily
give the Cox Theorem and the rules of probability without having
to modify the plausibility function.  In addition, merely adding
the axiom \ref{axadd}, then gives the full theory of expectation
for random variables as well as general unknown numbers.  In fact,
if we go to complex unknowns, with obvious complex versions of the
axioms, and assume that the unknowns form a $C^*-$algebra, as
specifying a $PV$ is equivalent to giving a state, it is known
that every state is a bounded linear map \cite{DIX}, so that the
usual analysis with measure theory follows from the representation
of bounded linear functions as integration with respect to a
finite measure.

Suppose that more generally we have a vector space $W$ and we are
interested in plausible values for members of a set $T$ of unknown
members of $W.$  Then, the obvious modification of the axioms
\ref{axhomog} and \ref{axcox} leads to the conclusion that if $X$
is in $T$ and $A,C$ belong to $E,$ then by (4) of corollary
(\ref{retractcor}) we find the obvious generalization of the Cox
Theorem again.  In fact, we can replace $W$ by a general module
$M$ over any possibly noncommutative ring $R,$ and with the
obvious modification of the axioms, we obtain the obvious
generalization of the Cox Theorem, where we simply replace
indicators by idempotents in ring $R$ and assume in addition that
$A$ and $C$ commute as idempotents in $R$ so as to make their
product again idempotent.  The main point here is that if $X$ is
an unknown member of $R$ and $v$ is in $W,$ then we must assume
that $PV(Xv|e)$ depends only on $PV(X|e)$ as a function of $X$
keeping $e$ and $v$ fixed.  As a consequence of this assumption we
find the general rule $$PV(Xv|e)=[PV(X|e)]v,$$ which combined with
(3) in Corollary (\ref{retractcor}) gives the general
multiplication rule:  $$PV(e_1Y|e_2)=PV(Y|e_1e_2)PV(e_1|e_2),$$ as
long as $e_1, e_2$ and $e_1e_2,$ are all idempotents, which is the
case if the two idempotents commute.  Here, $Y$ is an unknown
vector, so we must keep in mind that $PV(Y|e)$ is a member of $W$
whereas $PV(X|e)$ is a member of $R.$  Finally, if we take $T$ to
be a $C^*-$algebra and $R$ to be $C^*-$subalgebra and $P$ a
retraction of $T$ on $R,$ then Corollary (\ref{retractcor}) gives
us simple natural conditions for $P$ to be a conditional
expectation in $C^*-$algebra theory, that is conditions for $P$ to
be an $R-$ linear map.

Of course, we can produce examples of $PV$ functions by taking in
particular function algebras or even noncommutative
$C^*-$algebras. In particular, it is known that if $T$ is a
$C^*-$algebra with identity, and if we take for our set of
unknowns the set $S$ of self-adjoint members of $T,$ then any
state of the $C^*-$algebra restricted to $S$ will serve as a
consistent way of assigning plausible values which in fact satisfy
the general additivity of proposition (\ref{gensum}).  In fact, if
$T$ is any separable $C^*-$algebra, we can take the universal
representation and produce a state, $f,$ which will not vanish on
any nonzero positive element.  We then define the plausible value
$PV(X|A)=f(XA)/f(A),$ for any $X,A \in T$ such that $A$ is a
nonzero idempotent. In particular, if $T$ is commutative, then we
know that the states which are multiplicative are exactly the pure
states, which are the point evaluations under any representation
of such an algebra as an algebra of continuous functions on a
compact Hausdorff space. Thus, the assumption that plausible value
is generally additive is a reasonable assumption, whereas we see
that the additional assumption of multiplicativity would be too
restrictive. In general, it is known from Choquet theory that the
set of all states of a $C^*-$algebra is a compact convex subset of
the continuous dual of the algebra under the weak*-topology, and
that it is the closed convex hull of the pure states, as these
form the set of extreme points of that convex set \cite{DIX}.


\begin{thebibliography}{99}


\bibitem{Borel}E. Borel, A propos d'un trait\~e  de probabilit\'es, {\it Rev. Philos.}, {\bf 98} (1924), 321--336.


\bibitem{COX1}
R. T. Cox, Probability, frequency, and reasonable expectation,
\emph{ Am. J. Phys.}, \textbf{14}(1946), 1--13.

\bibitem{COX}
R. T. Cox, \emph{ The Algebra of Probable Inference}, The Johns
Hopkins Press, Baltimore, Maryland, 1961.


\bibitem{DIX}
J. DIXMIER, \emph{ Les $C^*-$Algebres et Leurs Representations},
Gauthier-Villars, Paris, France, 1969.


\bibitem{DUB}
D. Dubois and H. Prade, The logical view of conditioning and its
application to possibility and evidence theories.  \emph{
International Journal of Approximate Reasoning,} \textbf{4}(1),
23--46.

\bibitem{FINE}
T. L. Fine, \emph{ Theories of Probability-an Examination of
Foundations}, Academic Press, New York, New York, 1973.

\bibitem{FUCHS}
L. Fuchs, \emph{ Partially Ordered Algebraic Systems}, Pergamon
Press, Reading, Massachusetts, 1963.

\bibitem{HALP1}
J. Halpern,  A counterexample to theorems of Cox and Fine. \emph{
J. of Artificial Intelligence Research} \textbf{10}(1999),67--85.

\bibitem{HALP2}
J. Halpern, Technical Addendum, Cox's Theorem revisited, \emph{ J.
of of Artificial Intelligence Research}
\textbf{11}(1999),429--435.

\bibitem{HARDY}
M. Hardy, Scaled Boolean algebras, \emph{Adv. in Applied Math.},
\textbf{29}(2002), no. 2, 243--292.

\bibitem{HALMOS}
P. Halmos, \emph{ Lectures on Boolean Algebras}, van Nostrand,
Princeton, New Jersey, 1963.

\bibitem{KHH}
K. H. Hofmann and P. S. Mostert, \emph{ Elements of Compact
Semigroups}, C. E. Merrill Books, Columbus, Ohio, 1966.


\bibitem{JAYNES}
E. T. Jaynes, \emph{ Probability Theory-The Logic of Science},
Cambridge University Press, Cambridge, U.K., 2003.

\bibitem{Kolmogorov}A. N. Kolmogorov, {\it Foundations of the Theory of Probability}, (English translation of a 1933 German language original) Chelsea Publishing House, New York, 1950.

\bibitem{KOOPMAN}
B. O. Koopman, The axioms and algebra of intuitive probability,
\emph{ Annals of Mathematics}, \textbf{41}(1940), 269--292.

\bibitem{KRAFT}
C. H. Kraft, J. W. Pratt, and A. Seidenberg, Intuitive probability
on finite sets, \emph{Annals of Mathematical Statistics},
\textbf{30}(1959), 408-419.

\bibitem{Laplace}P. S. de Laplace, {\it Th\'eorie Analytique des Probabilit\'es} (2 volumes) Coucier Imprimeur, Paris, 1812.

\bibitem{LOOMIS}
L. H. Loomis, \emph{ The Lattice-theoretic Background of the
Dimension Theory of Operator Algebras}, Mem. A.M.S.,
\textbf{18}(1955).

\bibitem{Morgan}
A de Morgan, {\it Formal Logic: or the Calculus of Inference Necessary and Probable}, Taylor \& Watton, London 1847.

\bibitem{PARIS}
J. B. Paris, \emph{ The Uncertain Reasoner's Companion.}, Cambridge
University Press, Cambridge, U.K., 1994.

\bibitem{Poincare}
H. Poincar\'e, {\it Calcul de Probabilit\'es} (2nd edition),Gauthier-Villars, Paris, 1912.

\bibitem{Polya}G. P\'olya, {\it Mathematics and Probable Reasoning} Princeton University Press, Princeton, 1954.

\bibitem{SAVAGE}
L. J. Savage, \emph{ Foundations of Statistics}, Dover, New York,
New York, 1972.

\bibitem{SCOTT}
D. Scott, Measurement structures and linear inequalities,
\emph{Journal of Mathematical Psychology}, \textbf{1}(1964),
233--247.

\bibitem{SIVIA}
D. S. Sivia, \emph{ Data Analysis A Bayesian Tutorial}, Oxford
University Press, Oxford, U.K., 1996.

\bibitem{TIPLER}
F. J. Tipler, What about quantum mechanics? Bayes and the Born
interpretation, quant-ph/0611245.

\bibitem{TRIBUS1}
M. Tribus, \emph{ Rational Descriptions, Decisions and Designs},
Pergamon Press, New York, 1969.



\end{thebibliography}
\end{document}